\documentclass[12pt]{amsart}

\usepackage{amssymb}
\usepackage{mathrsfs}

\def\F{\mathbb{F}}

\def\su{{\subseteq}}
\def\la{{\langle}}
\def\ra{{\rangle}}

\def\w{{\omega}}

\def\dim{{\rm dim\,}}

\def\ad{\mbox{ad~}}

\def\Proof{\noindent{\sl Proof.}\ }
\def\qed{{\hfill $\Box$ \medbreak}}

\newtheorem{defi}{Definition}[section]
\newtheorem{thm}[defi]{Theorem}
\newtheorem{lem}[defi]{Lemma}
\newtheorem{cor}[defi]{Corollary}

\newtheorem{eg}[defi]{Example}

\begin{document}

\title[Engel condition on  enveloping algebras]{Engel condition on  enveloping algebras of Lie superalgebras}

\author{\textsc{Salvatore Siciliano}}
\address{Dipartimento di Matematica e Fisica ``Ennio De Giorgi", Universit\`{a} del Salento,
Via Provinciale Lecce--Arnesano, 73100--Lecce, Italy}
\email{salvatore.siciliano@unisalento.it}

\author{\textsc{Hamid Usefi}}
\address{Department of Mathematics and Statistics,
Memorial University of Newfoundland,
St. John's, NL,
Canada, 
A1C 5S7}
\email{usefi@mun.ca}

\thanks{The research of the second author was supported by NSERC of Canada}

\begin{abstract}
Let $L$  be a Lie superalgebra over a fled  of characteristic $p\neq 2$ with enveloping algebra $U(L)$ or let $L$ be a restricted Lie superalgebra over a field of characteristic $p > 2$ with restricted enveloping algebra $u(L)$. 
 In this note, we establish when $u(L)$ or $U(L)$ is bounded Lie Engel.
\end{abstract}

\date{\today}

  \maketitle

\section{Introduction}
Recall that an associative ring $R$ is said to satisfy the Engel condition if $R$ satisfies the identity
$$[x, \underbrace{y,\ldots, y}_n] = 0,$$
for some $n$. It follows from  Zel'manov's celebrated result about the restricted Burnside problem \cite{Z90}
that every finitely generated Lie ring satisfying the Engel condition 
is nilpotent. Kemer in \cite{kemer} proved that if $R$ is an associative algebra over a field of characteristic zero that satifies  the Engel condition then $R$ is Lie nilpotent. This result was later proved by Zel'manov in \cite{Z88}  for all Lie algebras. 
However these results fail in positive characteristic, see \cite{razm, Rips}. Nevertheless, Shalev in \cite{Sh} proved that every finitely generated associative algebra over a field of characteristic $p>0$ satisfying the Engel condition is Lie nilpotent. This result was further strengthened by Riley and Wilson in \cite{RW2} by proving that  if  $R$ is a $d$-generated associative $C$-algebra, where $C$ is a commutative ring,  satisfying the Engel condition of degree $n$, then $R$ is upper Lie nilpotent of class bounded by a function that depends only on $d$ and $n$. Hence, in the positive characteristic case one would need to assume that $R$ is also finitely generated.

Let $L=L_0\oplus L_1$ be a Lie superalgebra over a field $\F$ of characteristic $p\neq 2$
with bracket $(\,,)$.  The adjoint map of $x\in L$ is denoted  by  $\ad x$. We denote the enveloping algebra of $L$ by $U(L)$.
 In case $p=3$ we add the condition $((y,y),y)=0$, for every $y\in L_1$.  This identity is necessary to embed $L$ in $U(L)$.
 
 The Lie bracket of $U(L)$ is denoted by $[a,b]=ab-ba$, for every $a, b\in U(L)$.
We are interested to know when $U(L)$ satisfies the Engel condition. Note that the Engel condition is a non-matrix identity, that is a polynomial identity not satisfied by the algebra $M_2(\F)$ of $2\times 2$ matrices over $\F$. The conditions for which $U(L)$ satisfies a non-matrix identity are given in \cite{BRU}. 
It follows from Zel'manov's Theorem \cite{Z88}  that over a field of characteristic zero $U(L)$ satisfies the Engel condition if and only if $U(L)$ is Lie nilpotent. The characterization of $L$ when $U(L)$ is Lie nilpotent over any field of characteristic not 2 is given in  \cite{BRU}. Hence, we have
\begin{cor}\label{UL-Engel-zero-char}
Let $L=L_0\oplus L_1$ be a Lie superalgebra over a field of characteristic zero. The following conditions are equivalent:
\begin{enumerate}
\item $U(L)$ is Lie nilpotent;
\item $U(L)$ is bounded Lie Engel;
\item $L_0$ is abelian, $L$ is nilpotent, $(L,L)$ is finite-dimensional,  and  either   $(L_1, L_1)=0$ or $\dim\, L_1\leq 1$ and $(L_0, L_1)=0$.
\end{enumerate}
\end{cor}

However this result is no longer true in positive characteristic as our following theorem shows (see also Example \ref{counterexample}).

\begin{thm}\label{UL-Engel-positive-char}
Let $L=L_0\oplus L_1$ be a Lie superalgebra over a field of characteristic $p\geq 3$. The following conditions are equivalent:
\begin{enumerate}
\item $U(L)$ is bounded Lie Engel;
\item $U(L)$ is PI,    $L_0$ is abelian,     $\ad x$ is  nilpotent  for every $x\in L_0$, and either   $(L_1, L_1)=0$ or $\dim\, L_1\leq 1$ and $(L_0, L_1)=0$;
\item $U(L)$ is PI, $L_0$ is abelian, $L$ is nilpotent,   and  either   $(L_1, L_1)=0$ or $\dim\, L_1\leq 1$ and $(L_0, L_1)=0$.
\end{enumerate}
\end{thm}

Note that the above theorem does not follow from Zel'manov or  Riley and Wilson's results because $U(L)$ is not necessarily finitely generated. 

Now let $L=L_0\oplus L_1$ be a restricted Lie superalgebra over a field  of characteristic $p> 2$ with enveloping algebra $u(L)$. In our next result we characterize $L$ for which $u(L)$ satisfies the Engel condition. 
Our results complement the  results of \cite{U1, U2} where it is determined when  $u(L)$ satisfies a 
non-matrix identity or when $u(L)$ is Lie solvable, Lie nilpotent, or Lie super-nilpotent. Similar results for group rings and enveloping algebras of restricted Lie algebras were carried out in \cite{BG,PPS} and  \cite{RS}, respectively.

\begin{thm}\label{uL-Lie Engel}
Let $L=L_0\oplus L_1$ be a restricted Lie superalgebra over a  field of characteristic $p>2$.
The following conditions are equivalent:
\begin{enumerate}
\item  $u(L)$ is bounded Lie Engel;
\item $u(L)$ is PI,    $(L_0, L_0)$ is $p$-nilpotent,  there exists an integer $n$ such that $(\ad x)^n=0$ for every $x\in L_0$, and either $(L_1, L_1)$ is $p$-nilpotent or $\dim L_1\leq 1$  and $(L_1, L_0)=0$;
\item  $u(L)$  is PI, $L$ is  nilpotent,   $(L_0, L_0)$ is $p$-nilpotent,
and either $(L_1, L_1)$ is $p$-nilpotent or $\dim L_1\leq 1$  and $(L_1, L_0)=0$.
\end{enumerate}
\end{thm}

We refer the reader to \cite{B92} for basic background about Lie superalgebras and their enveloping algebras and to \cite{BRU, U1, U2} for notation. 

\section{Proof of Theorem \ref{UL-Engel-positive-char}}\label{ordinary}

Throughout the paper, all Lie superalgebras are defined over a field of characteristic $p\geq 3$.

Let $\mathfrak{D}$ be a subset of an associative algebra $\mathfrak{U}$ over a field $\F$. Recall that $\mathfrak{D}$ is called \emph{weakly closed} if for every pair of elements $(a, b)\in \mathfrak{D}\times \mathfrak{D}$, there exists an element $\gamma(a, b)\in \F$ such that $ab+\gamma(a, b)ba\in \mathfrak{D}$. This notion is applicable to our setting with a Lie superalgebra $L$. We take $\mathfrak{U}$ to be the associative subalgebra of $\mbox{End}_{\F}(L)$ generated by all $\ad x$, where $x\in L$, and $\mathfrak{D}$ to be the subset of $\mathfrak{U}$ consisting of all $\ad x$, where $x$ is a homogeneous element in $L$. Then $\mathfrak{D}$ is weakly closed. Now we recall the Jacobson's Theorem on weakly closed sets.

\begin{thm}[Jacobson {\cite{J79}}]\label{jacobson}
Let $\mathfrak{D}$ be a weakly closed subset  of an  associative algebra $\mathfrak{U}$ 
of linear transformations  of a finite-dimensional vector space $V$ over $\F$. Assume that every $T\in \mathfrak{D}$ is  nilpotent. Then the non-unital associative subalgebra generated by $\mathfrak{D}$ is associative nilpotent.
\end{thm}

Next we recall Corollary 2.5 from \cite{P92}. Note that the derived subalgebra of $L$ will be also denoted by $L'$.

\begin{thm}\label{Petrogradsky-ordinary}
Let $L=L_0\oplus L_1$ be a Lie superalgebra over a field of  characteristic $p>2$.
Then $U(L)$ is PI if and only if there exist homogeneous ideals $B\su A\su L$ such that
\begin{enumerate}
\item $\dim L/A<\infty, \dim B< \infty$;
\item $A'\su B$;
\item $B=B_1$; 
\item All inner derivations $\ad z, z\in L_0$, defined over $L$ are algebraic and their degrees
 are bounded by some constant.
 \end{enumerate}
\end{thm}

\begin{lem}[\cite{BRU}] \label{L0-central}
 If $U(L)$ satisfies a non-matrix polynomial identity then $(L_0, L_0)=(L_0, L_1, L_1)=0$.
\end{lem}

\begin{lem} \label{L1-L1-ordinary}
If  $U(L)$ is bounded Lie Engel then either $(L_1, L_1)=0$ or $\dim L_1\leq 1$ and $(L_0, L_1)=0$.
\end{lem}
\Proof The proof follows exactly as  in Lemma 5.1 of \cite{BRU}.\qed

Note that  unlike the characteristic zero case,  $U(L)$ satisfying the Engel condition does not necessarily imply that $U(L)$ is Lie nilpotent. In fact, in  characteristic $p\geq 3$, $U(L)$ can be $p$-Engel and yet $(L, L)$ be infinite-dimensional.

\begin{eg}\label{counterexample}
\emph{ Let $L=L_0\oplus L_1$, where $(L_0, L_0)=(L_1, L_1)=0$, $L_0$ has a basis $x_1, x_2, \cdots$ and
$L_1$ has a basis $y, z_1, z_2, \cdots$ with $z_i=(x_i, y)$. Let $m$ be a positive integer and set 
 $v=x_2\cdots x_m$. Note that $[x_i, y, y, y]=0$, for all $x_i$. Thus, by the
Leibniz formula we have,
\begin{align*}
[x_1\cdots x_m,_p y]&=x_1[v,_p y]+p[x_1, y][v,_{p-1} y]+ {p \choose 2}[x_1, y, y][v,_{p-2} y]\\
&=x_1[v,_p y]=\cdots=x_1\cdots x_{m-1}[x_m,_p y]=0.
\end{align*} 
It follows that $U(L)$ is $p$-Engel.}
\end{eg}

We just recall the following identity that follows from super-Jacobi identity:

\begin{align}\label{Jacobi}
(\text{ad } z)^2=\frac{1}{2} \text{ad } (z, z), \quad \text{ for every } z\in L_1.
\end{align}

\begin{lem}\label{L-nilpotent-ordinary}
Suppose  $U(L)$ is PI,    $(L_1, L_1)=(L_0, L_0)=0$, and  $\ad x$ is a nilpotent transformation on $L$, for every $x\in L_0$.  Then $L$ is nilpotent.
 \end{lem}
 \Proof 
 Note that, by Theorem \ref{Petrogradsky-ordinary},    $L$ contains a homogeneous ideal $A$ of finite codimension such that $A'$ is finite-dimensional. First we show that $A$ is nilpotent. 
Note that the restriction of $\ad x$  to  $A'$  is a nilpotent transformation acting on a finite-dimensional vector space, for every $x\in A_0$.
It follows from Theorem \ref{jacobson} that  the (non-unital) associative algebra generated by all $\ad x$ with $x\in A_0$ acting on $A'$  is associative nilpotent. This means that
 $$(A, A, \underbrace{A_0, \ldots, A_0}_t)=0,$$ 
 for some $t$. Since $(L_1, L_1)=0$, we deduce that $A$ is nilpotent.
To prove $L$ is nilpotent, we argue by induction on $\dim L/A$. 
By Lemma 5.4 in  \cite{SU11}, it suffices  to show that $L/A'$ is nilpotent. So we can replace $L$ with $L/A'$ and assume that $A$ is an abelian ideal of $L$ of finite codimension.
Let $z$ be a homogeneous element in  $L\setminus A$ and denote by  $N$ the ideal of $L$
generated by $A$ and $z$. If $z\in L_1$ then we deduce from Equation \eqref{Jacobi} and the hypothesis
$(L_1, L_1)=0$ that $(A, z, z)=0$. It is now easy to see that $N$ is nilpotent. On the other hand, if $z\in L_0$ then 
 $\ad z$ is a nilpotent transformation on $L$.  Since $A$ is abelian, we can observe that $N$ is nilpotent. 
Now we note that  $\dim\, L/N<\dim\, L/A$ and it follows by induction that $L$ is nilpotent.
  \qed

\noindent\emph{Proof of Theorem \ref{UL-Engel-positive-char}.} 
The implication $(1)\Rightarrow (2)$  follows from Lemmas \ref{L0-central} and    \ref{L1-L1-ordinary} while 
$(2)\Rightarrow (3)$ follows from Lemma \ref{L-nilpotent-ordinary}.
It remains to prove $(3)\Rightarrow (1)$. Note that, by Theorem 1.1 of \cite{BRU}, $U(L)$ satisfies a non-matrix PI and thus $R=[U(L), U(L)]U(L)$ is nil of bounded index, say $p^m$. Since $L$ is  nilpotent, there exists an integer $n$ such that $L_0^{p^n}$ is contained in the centre $Z(L)$ of $L$. Let $X$ be an ordered basis for $L_0$ and $Y$ an ordered basis for $L_1$.
Let $w$ be an element in the augmentation ideal $\w(L)$ of $U(L)$. Then by the PBW Theorem, $w$ is a linear combination of PBW monomials of the form
$x_1^{a_1}\cdots x_i^{a_i}y_1\cdots y_j$. 

First  suppose that $(L_1, L_1)=0$. Then, modulo $R$, $w^{p^n}$ is a linear combination of monomials of the form 
\begin{align*}
x_1^{a_1p^n}\ldots x_i^{a_ip^n}.
\end{align*}
Since $L_0^{p^n}\su Z(L)$, we deduce that $w^{p^n}=u+v$, where $u$ is a central element in $U(L)$ and $v\in R$.  Hence,  $w^{p^{m+n}}=u^{p^m}$ is a central element in $U(L)$. Clearly, $m$ and $n$ are independent of $w$ and so  $U(L)$ is $p^{m+n}$-Engel in this case.

The case when $\dim L_1\leq 1$ and $(L_0, L_1)=0$ can be handled similarly.
\qed

\section{Proof of Theorem \ref{uL-Lie Engel}}\label{restricted}

We recall that Engel's Theorem holds for Lie superalgebras (see \cite{sch}, for example).
\begin{thm}[Engel's Theorem]\label{Engel's Theorem}
Let $L$ be a finite-dimensional Lie superalgebra such that
$\text{ad } x$ is nilpotent, for every homogeneous element $x\in L$.
Then $L$ is nilpotent.
\end{thm}

\begin{lem} \label{L1-L1}
If  $u(L)$ is bounded Lie Engel then either $(L_1, L_1)$ is $p$-nilpotent  or $\dim L_1\leq 1$ and 
$(L_0, L_1)=0$.
\end{lem}
\Proof The proof follows exactly as  in Lemma 4.1 of \cite{U2}.\qed

\begin{thm}[\cite{P92}]\label{Petrogradsky-rest}
Let $L=L_0\oplus L_1$ be a restricted Lie superalgebra.
Then $u(L)$ satisfies a  PI if and only if there exist homogeneous restricted  ideals $B\su A\su L$ such that
\begin{enumerate}
\item $L/A$ and $B$ are both finite-dimensional.
\item $A'\su B$, $B'=0$.
\item The restricted Lie subalgebra  $B_0$ is $p$-nilpotent.
 \end{enumerate}
\end{thm}

\begin{thm}[\cite{U1}]\label{uL-non-matrix}
Let $L=L_0\oplus L_1$ be a restricted Lie superalgebra over a perfect field and denote by $M$  the  subspace spanned by all $y\in L_1$ such that $(y,y)$ is $p$-nilpotent.
The following statements are equivalent:
\begin{enumerate}
\item $u(L)$ satisfies a non-matrix PI.
\item $u(L)$ satisfies a PI, $(L_0, L_0)$ is $p$-nilpotent,  $\dim L_1/M\leq 1$,   $(M, L_1)$ is $p$-nilpotent, and  $(L_1, L_0)\su M$.
\item The commutator ideal of $u(L)$ is nil of bounded index.
\end{enumerate}
\end{thm}

\begin{lem}\label{L-nilpotent}
Suppose that $u(L)$ is PI,    $(L_1, L_1)$ is $p$-nilpotent or $(L_0, L_1)=0$,  and  $\ad x$ is a nilpotent transformation on $L$, for every $x\in L_0$.  Then $L$ is nilpotent.
 \end{lem}
 \Proof   Note that, by Theorem   \ref{Petrogradsky-rest},    $L$ contains a homogeneous ideal $A$ of finite codimension such that $A'$ is finite-dimensional and $(A', A')=0$. First we show that $A$ is nilpotent. Let $a\in A'$  and $y\in A_1$. By Equation \eqref{Jacobi}, we have
  $$
  (a, y, y)=0.
  $$
  Consider now the set $\mathfrak{D}$ of all linear transformations $\ad z$ acting on  $A'$, where $z$ is a homogeneous element in $A$.  Since $\mathfrak{D}$ is  weakly closed and consisting of nilpotent linear transformations, we deduce, by Theorem \ref{jacobson}, that the non-unital associative algebra generated by $\mathfrak{D}$ is associative nilpotent. This means that   
  $$
  (A', \underbrace{A, \ldots, A}_t)=0,
  $$
   for some $t$. Hence, $A$ is nilpotent.
To prove that $L$ is nilpotent, we use induction on $\dim L/A$.
By Lemma 5.4 in  \cite{SU11}, it suffices  to show that $L/\la A'\ra_p$ is nilpotent. 
So we replace $L$ with $L/\la A'\ra_p$ and assume that $A$ is an abelian ideal of $L$ of finite codimension.
We claim that $H=L/A$ is nilpotent.  Note that  either $(H_0, H_1)=0$ or 
 $(H_1, H_1)$ is $p$-nilpotent. If $(H_0, H_1)=0$ then clearly every $\text{ad } x$ with $x\in H_1$ is nilpotent. On the other hand, if  $(H_1, H_1)$ is $p$-nilpotent then, by Equation \eqref{Jacobi}, every $\text{ad } x$ with $x\in H_1$ is nilpotent. Hence,  Theorem \ref{Engel's Theorem} applies and we deduce that $H$ is nilpotent. This means that
$\gamma_{c+1}(L)\su A$, where $c$ is the nilpotency class of $H$. 
Let $z$ be a homogeneous element in  $\gamma_c(L)$ and denote by  $N$ the ideal of $L$
generated by $A$ and $z$. If $z\in L_1$ then, by Equation \eqref{Jacobi}, $(A, z, z)=0$. 
On the other hand, $\ad z$ is nilpotent if $z\in L_0$. We deduce that $N$ is nilpotent. 
Now, $\dim\, L/N<\dim\, L/A$ and it follows by induction that $L$ is nilpotent.
\qed

\noindent\emph{Proof of Theorem \ref{uL-Lie Engel}.}  
The implication $(1)\Rightarrow (2)$  follows from Lemma \ref{L1-L1} and Theorem \ref{uL-non-matrix} while 
$(2)\Rightarrow (3)$ follows from Lemma \ref{L-nilpotent}.
It remains to prove $(3)\Rightarrow (1)$. Note that, by Theorem \ref{uL-non-matrix}, $u(L)$ satisfies a non-matrix PI and thus $R=[u(L), u(L)]u(L)$ is nil of bounded index. Since $L$ is  nilpotent, there exists an integer $s$ such that $L_0^{p^s}\su Z(L)$. Let $X$ be an ordered basis for $L_0$ and $Y$ an ordered basis for $L_1$.
Let $w\in \w(L)$. Then, by the PBW Theorem, $w$ is a linear combination of PBW monomials of the form
$x_1^{a_1}\cdots x_i^{a_i}y_1\cdots y_j$. Let $z_j=\frac{1}{2}(y_j, y_j)$.
Note that $y_j^p= z_j^{\frac{p-1}{2}} y_j$. 

Now suppose that $(L_1, L_1)$ is $p$-nilpotent. Then there exists an integer $n\geq s$ 
such that $y_j^{p^n}=0$, for all $y_j\in Y$.   Now  note that, module $R$, $w^{p^n}$ is a linear combination of monomials of the form 
\begin{align*}
x_1^{a_1p^n}\ldots x_i^{a_ip^n}.
\end{align*}
Since, $L_0^{p^n}\su Z(L)$, we deduce that $w^{p^n}\in Z(u(L))+R$. Furthermore,  $R$ is nil of bounded index, say $p^m$. Hence,  $w^{p^{m+n}}\in Z(u(L))$. Clearly, $m$ and $n$ are independent of $w$ and so  $u(L)$ is $p^{m+n}$-Engel in this case.
 
 On the other hand, if $(L_1, L_1)$ is not  $p$-nilpotent then, by the hypothesis, we must have $\dim L_1= 1$ and $(L_1, L_0)=0$.  Suppose that $L_1$ is spanned by $y\in L_1$ and let $u$ be a PBW monomial of the form $x_1^{a_1}\cdots x_i^{a_i}y$. We have 
 $u^{p^s}=x_1^{a_1p^s}\ldots x_i^{a_ip^s}y^{p^s}\in Z(u(L))+R$. Similarly, if $v$ is a PBW monomial of the form $x_1^{a_1}\cdots x_i^{a_i}$ then $v^{p^s}\in Z(u(L))+ R$. We deduce that $w^{p^s}\in Z(u(L))+ R$. Clearly, $m$ and $s$ are independent of $w$ and so  $u(L)$ is $p^{m+s}$-Engel, completing the proof.  \qed

\end{document}